\documentclass[12pt]{amsart}
\usepackage[english,francais]{babel}
\usepackage[T1]{fontenc}
\usepackage{amsmath,amsfonts}
\usepackage{amsthm}
\usepackage{pifont}
\usepackage[all]{xy}
\usepackage{tikz}
\usepackage{extarrows}

\usepackage{geometry}
\usetikzlibrary{fit,calc,positioning,decorations.pathreplacing,matrix}

\newtheorem{theo}{Theorem}[section]

\newtheorem{Rq}{Remark}

\newtheorem{ex}{Example}
\newtheorem{cor}{Corollary}

\author{K. Boutahir and, Y. Rami}
\address{D\'epartement de Math\'ematiques \& Informatique,\\
Universit\'e  My Ismail, B. P. 11 201 Zitoune, Mekn\`es, Morocco,}
\email{khalid.boutahir@edu.umi.ac.ma and yousfoumadan@gmail.com}

\title{On L.S.-category of a family of rational elliptic spaces\\
     }
\date{11/10/2013}


\begin{document}
\maketitle

\selectlanguage{english}
\begin{abstract}
Let X be a finite type simply connected rationally elliptic CW-complex
with Sullivan minimal model $(\Lambda V, d)$ and let $k \geq 2$ the biggest integer
such that $d=\sum \limits_{\underset{}{i\geq k}}d_i$ with $d_i(V ) \subseteq \Lambda^iV$.

In \cite{murillo02} the authors showed that if $(\Lambda V,d_k)$ is morever
elliptic then  $cat(\Lambda V,d)=(k-2)dimV^{even} + dimV^{odd}.$

Our work focuses on the estimation of L.S.-category of such spaces in
the case when $k=3$ and when $(\Lambda V,d_3)$ is not necessarily
elliptic.
\end{abstract}


\section{Introduction \label{intro}}

Let X be a finite type simply connected  CW-
complex with Sullivan minimal model $(\Lambda V, d)$ and let $k \geq 2$ the biggest integer
such that $d=\sum \limits_{\underset{}{i\geq k}}d_i$ with $d_i(V ) \subseteq \Lambda^iV$
 and $dim(V ) < \infty$.\\

Consider on $(\Lambda V, d)$  the filtration given by $$F^p = \Lambda^{\geq (k-1)p}V =\bigoplus_{i=(k-1)p}^{\infty}\Lambda^i V.$$
$F^p$ is preserved by the differential d and satisfies $F^p(\Lambda V)\otimes F^q(\Lambda V)\subseteq F^{p+q}(\Lambda V)$, $\forall p, q \geq 0$, so it is a filtration  of differential graded algebras. Also,
since $F^0=\Lambda V$ and $F^{p+1}\subseteq F^p$ this filtration is decreasing and bounded, so it induces a convergent spectral sequence.
Its $0^{th}$-term is
$$E_0^{p,q}=\bigg(\frac{F^p}{F^{p+1}}\bigg)^{p+q} = \bigg(\frac{\Lambda^{\geq (k-1)p}V}{\Lambda^{\geq (k-1)(p+1)}V}\bigg)^{p+q}.$$
Hence, we have the identification:
\begin{equation} \label{ssk}
E_0^{p,q}=
\big(\Lambda^{p(k-1)}V\oplus\Lambda^{p(k-1)+1}V\oplus...\oplus\Lambda^{p(k-1)+k-2}V
\big)^{p+q}\,\,\,\,\,\,\,\,\,
\end{equation}

In this general situation, the $1^{th}$-term is the graded algebra $\Lambda V$ proveded with a differential $\delta $, which is'nt necessarely a derivation on the set $V$ of generators (see $\S3$). That is $(\Lambda V, \delta )$ is  a commutative differential graded algebra, but it is not a Sullivan algebra. The spectral sequence is therefore:
$$H^{p,q}(\Lambda V, \delta)\Rightarrow H^{p+q}(\Lambda V, d).$$
Hence if $dim(V)< \infty$ and $(\Lambda V,\delta)$ has finite dimensional cohomology, then  $(\Lambda V,d)$ is elliptic. This gives a new family
of rationally elliptic spaces for which $d=\sum \limits_{\underset{}{i\geq k}}d_i$ 

Recall first that in \cite{murillo02} the authors gives the  explicit formula $cat(\Lambda V,d) = dimV^{odd} + (k-1)dimV^{even}$ of L.-S. category for a  minimal Sullivan model $(\Lambda V,d)$ satisfying the restrictive condition : $(\Lambda V,d_k)$ is also elleptic. 

It is important to note also that their algorithm that induces the fundamental class of $(\Lambda V,d)$ from that of $(\Lambda V,d_k)$ corresponds to the progress of a cocycle that survive to term $E_{\infty}$ (cf. Remark 1). 

The main result of this work is a project of  determination of an explicit formula for  $cat(\Lambda V,d)$ with $(\Lambda V,d)$ being elliptic
and $(\Lambda V,d_k)$ not elliptic, completing the formula given by
L. Lechuga and A. Murillo in \cite{murillo02}.

In what follow, we consider the case where $d=\sum \limits_{\underset{}{i\geq 3}}d_i$, that is where $k=3$ and $N$ designate  the formal dimention of $(\Lambda V, d)$. With the notation as above, our first result reads:

\begin{theo} If  $(\Lambda V,d)$ is  elleptic and   $H^N(\Lambda V, \delta ) = \mathbb{Q}.\alpha$ is one dimentional, 
then
$ cat_0(X)=cat(\Lambda V,d)= sup\{ k\geq 0, \; \alpha = [\omega _0]  \; with \;  \omega _0 \in \Lambda^{\geq k} V\}.$
\end{theo}
\medskip

Let  $(\Lambda W,d)$ a minimal Sullivan model of $(\Lambda V,\delta )$. If $dim(W)<\infty $ 
then (\cite{murillo94}) $(\Lambda W,d)$  is a Gorenstein algebra  and so is $(\Lambda V,\delta )$.  If additionaly $dim H(\Lambda V,\delta )<\infty$, then (\cite{felix82}) $(\Lambda W,d)$ is elliptic and so  its L.S. category is finite. It follows (\cite[Th. 29.15]{felix01}) that $Mcat(\Lambda V, \delta )<\infty $. Hence (\cite[Th. 3.6]{felix88}) $H(\Lambda V, \delta )$ is a Poincaré Duality algebra. There follow the
\begin{cor}
Let  $(\Lambda W,d)$ a minimal Sullivan model of $(\Lambda V,\delta )$. If $dim(W)<\infty $ and $dimH(\Lambda V, \delta )<\infty$ then
$cat_0(X)= sup\{ k\geq 0, \; \alpha = [\omega _0]  \; with \;  \omega _0 \in \Lambda^{\geq k} V\}.$
\end{cor}

\begin{Rq}

Now if $dimH^N(\Lambda V, \delta ) > 1$,   the technique used to show Theorem 1.1 can be adapted to have a similar result under this general hypothesis. The procedure is as follows:

Note first that in the proof of Theorem 1.1,  the algorithm  applied to the representative $\omega _0$ of the generating class of $H^N(\Lambda V,\delta)$ resulted in one of the fundamental class of  $H(\Lambda V, d)$ because $\omega _0$ is a cocycle which survives to $E_{\infty}$ in the spectral sequence.

On the other hand, since $dim(V)<\infty $,  we have   $dimH^N(\Lambda V, \delta )<\infty $, with $N$ being the formal dimension of $(\Lambda V,d)$. Since the  filtration induces on cohomology a graduation such that $H^N(\Lambda V,\delta) = \oplus _{p+q=N}H^{p,q}(\Lambda V,\delta)$, there is a basis $\{\alpha_1,...,\alpha_m\}$ of $H^N(\Lambda V,\delta)$ with $\alpha _i\in H^{p_i,q_i}(\Lambda V,\delta)$, $(1\leq i \leq m)$. That is, $\alpha _i = [(\omega _0^i,\omega _1^i)], \; \hbox{where}\;  (\omega _0^i,\omega _1^i)\in \Lambda ^{2p_i}V \oplus \Lambda ^{2p_i+1}V$.
Also since $(\Lambda V,d)$ is elleptic, there exist a unique  $j$ such that some $\alpha _j\in H^{p_j,q_j}(\Lambda V,\delta)$   survives to $E_{\infty}$ and consequently induces an representative of the fundamental class of $(\Lambda V,d)$. Explicitly, the corresponding obstructions $[a_2^0]=0$, $[a_3^1]=  0$, $\ldots $,  $[a_{t_j +l_j }^{t_j +l_j -2}] =0$ are necessary satisfied.

Now,  applying to   $(\omega _0^j,\omega _1^j)\in \Lambda ^{2p_j}V \oplus \Lambda ^{2p_j+1}V$,  the same role as  that applied to  $ \omega_0$ in the case of the first inequality (see \S 4) we obtain   an $\omega _{l_j + t_j -1}\in \Lambda ^{\geq 2p_j}V$ (resp. $\omega _{l_j + t_j -1}\in \Lambda ^{\geq 2p_j +1}V$) if $\omega _0^j\not = 0$  (resp.  if $\omega _0^j = 0$) representing the fundamental class of $(\Lambda V,d)$. It follows  that $e_0 (\Lambda V,d)\geq 2p_j$ (resp. $e_0 (\Lambda V,d)\geq 2p_j+1$). 

For the other inequality, any representative $ \omega\in\Lambda^{\geq s}V$ (where $s=e_0(\Lambda V,d)$) of the fundamental class of $(\Lambda V,d)$ induces by the same way, a representative  $(\omega_0, \omega_1)\in \Lambda ^{\geq s}V$ of a certain non zero class in $H^N(\Lambda V, \delta) = \oplus _{p+q=N}H^{p,q}(\Lambda V,\delta)$. By convergence of the spectral sequence $[\omega]$ correspond to a basis element of $E_{\infty}^{s,N-s}$ whch is one-dimentional, by ellepticity. It follow that in $E_{2}^{s,N-s}$ there is an element which survives to $E_{\infty}^{s,N-s}$. Hence with notations above, $s=2p_j$ or $s=2p_j +1$.  Therefore $e_0 (\Lambda V,d) = 2p_j$ or $e_0 (\Lambda V,d) = 2p_j +1$.
\end{Rq} 

With the notation of the previous remark we can therefore state the following generalization of the previous theorem
\begin{theo}
 If $(\Lambda V,d)$ is elleptic and $dimH^N(\Lambda V,d)=m$ with basis $\{\alpha _1, \ldots , \alpha _m\}$. Then 
$cat_0(X)= cat(\Lambda V,d) =r_j$ with $r_j=2p_j$ or else $r_j=2p_j +1$.
\end{theo}

\section{Basic facts and properties}

Let $\mathbb{K}$ be a field of characteristic $\neq$ 2.

  A Sullivan algebra is a free commutative differential graded algebra
(cdga for short) ($\Lambda V$, d) (where $\Lambda V$
=Exterior($V^{odd}$)$\otimes$ Symmetric($V^{even}$)) generated by
the graded $\mathbb{K}$-vector space
$V=\bigoplus_{i=0}^{i=\infty}V^i$ which has a well ordered basis
$\{x_\alpha\}$ such that $dx_\alpha$ $\in$ $\Lambda V_{<\alpha}$.
Such algebra is said minimal if $deg(x_\alpha)< deg(x_\beta)$
implies $\alpha <\beta$. If $V^0 = V^1$ = 0 this is equivalent to
saying that $d(V )\subseteq \bigoplus_{i=2}^{i=\infty}\Lambda^i
V$.\\
\newline
 A Sullivan model for a commutative differential graded algebra (A,
d) is a quasi- isomorphism (morphism inducing isomorphism in
cohomology) $(\Lambda V, d)\longrightarrow (A, d)$ with source, a
Sullivan algebra. If $H^0(A) = K$, $H^1(A) = 0$ and $dim(H^i (A, d))
<\infty$ for all $i\geq 0$, then \cite[Th.7.1]{halperin92}, this minimal
model exists. If X is a topological space any (minimal) model of the
algebra $C^*(X,\mathbb{K})$ is said a Sullivan (minimal) model of
X.\\
\newline
 The differential $d$ of any
element of V is a "polynomial" in $\Lambda V$ with no linear term. A
model ($\Lambda V$,$d$) is $elliptic$ if both V and $H^*$($\Lambda
V$,$d$) are finite dimentional spaces (see for example \cite{felix01}) .\\

For an elliptic space with model ($\Lambda V$, d) the formal
dimension $N$, i.e., the largest $n$ for which $H^n(\Lambda V,
d)\neq 0$, is given by \cite{halperin77}
$$N = dim V^{even} - \sum _{i =1}^{dim V}(-1)^{|x_i|}|x_i|$$
An element $0\neq \omega \in H^N(\Lambda V, d)$
is called a fundamental or top class of $(\Lambda V, d)$.\\

In \cite{halperin92} S.Halperin associated to any minimal model ($\Lambda V$, d) a pure model $(\Lambda V, d_{\sigma})$ defined as follows:\\
If $Q=V^{even}$ and $P=V^{odd}$ then
\begin{center}
$(\Lambda V, d_{\sigma})=(\Lambda Q\otimes\Lambda P, d_{\sigma})$; \,\,\,\,\,\,\, $d_{\sigma}(Q)=0 $ \,\,\,\,\,and \,\,\,\,\, $(d-d_{\sigma})(P)\subseteq\Lambda Q\otimes\Lambda^+P$

\end{center}

This model is related to $(\Lambda V, d)$ via the odd spectral sequence
$$H^{p,q}(\Lambda V, d_{\sigma})\Rightarrow H^{p+q}(\Lambda V, d)$$




The main result using this algebra and due to S. Halperin (\cite{halperin77}) shows that in the rational case, if $dim(V)<\infty$, then:
$$dim(H(\Lambda V,d))<\infty \Leftrightarrow dim(H(\Lambda V,d_{\sigma}))<\infty$$

 If X is a topological space, cat(X) is the least integer n such that X is covered
by n+1 open subset $U_i$, each contractible in X. It is an invariant of homotopy
type (c.f. \cite{felix01}). In \cite{felix}  Y. Félix,
S. Halperin and J.M. Lemaire showed that for Poincaré duality spaces, the
rational LS-category coincide with the rational Toomer invariant denoted $e_0(X)$.\\

By \cite[Lemma 10.1]{felix82} the Toomer invariant of a minimal
model $e_0(\Lambda V, d)$ is the
 largest integer s for which there is a non trivial cohomology class
 in $H^*(\Lambda V, d)$ represented by a cycle in $\Lambda^{\geq s}
 V$. As usual, $\Lambda^sV$ denotes the elements in $\Lambda V$ of
 $"wordlength" s$. For more details \cite{felix01},  \cite{halperin83}, \cite{sullivan78} are standard
 references.\\

 In \cite{murillo93} A. Murillo gave an expression of the fondamental class of $H(\Lambda V,d)$ in the case where $(\Lambda V,d)$ is a pure model. We recall it here:\\
 Assume $dim V <\infty$, choose
homogeneous basis $\{x_1,...,x_n\}$, $\{y_1,...,y_m\}$ of $V^{even}$
and $V^{odd}$ respectively, and write
 $$dy_j = a^1_j x_1 + a^2_j x_2+...+a^{n-1}_j x_{n-1}+a^n_j x_n\,\,\,\,\,\, j=1,2,...m,$$
where each $a^i_j$ is a polynomial in the variables
$x_i,x_{i+1},...,x_n$, and consider the matrix,

$$A=\begin{pmatrix}
\begin{tikzpicture}
\node (a) at (0,0) {$a_1^1$}; \node (b) at (0.5,0) {$a_1^2$}; \node
(c) at (1.5,0) {$a_1^n$};

 \node (d) at (0,-0.5) {$a_2^1$}; \node (e) at
(0.5,-0.5) {$a_2^2$}; \node (f) at (1.5,-0.5) {$a_2^n$};

 \node (g) at
(0,-1.5) {$a_m^1$}; \node (h) at (0.5,-1.5) {$a_m^2$}; \node (i) at
(1.5,-1.5) {$a_m^n$};

         \draw[dotted] (b.east) -- (c.west);
         \draw[dotted] (d.south) -- (g.north);
         \draw[dotted] (e.east) -- (f.west);
         \draw[dotted] (e.south) -- (h.north);
         \draw[dotted] (h.east) -- (i.west);
         \draw[dotted] (f.south) -- (i.north);
\end{tikzpicture}
\end{pmatrix}$$

 For any $1 \leq j_1<...<j_n\leq m$, denote by
$P_{j_1...j_n}$ the determinant of the matrix of order n formed by
the columns $i_1, i_2, ..., i_n$ of A:

$$\begin{pmatrix}
\begin{tikzpicture}
\node (a) at (0,0) {$a_{j_1}^1$}; \node (b) at (1.5,0)
{$a_{j_1}^n$}; \node (c) at (0,-1.5) {$a_{j_n}^1$}; \node (d) at
(1.5,-1.5) {$a_{j_n}^n$};
         \draw[dotted] (a.east) -- (b.west);
         \draw[dotted] (b.south) -- (d.north);
         \draw[dotted] (c.east) -- (d.west);
         \draw[dotted] (a.south) -- (c.north);
         \draw[dotted] (a.south east) -- (d.north west);
\end{tikzpicture}
\end{pmatrix}$$

Then (see \cite{murillo93}) if $dim H^*(\Lambda V,d)<\infty$ the
element $\omega\in\Lambda V$

\begin{equation} \label{omega}
\omega=\sum\limits_{\underset{}{1 \leq j_1<...<j_n\leq
m}}(-1)^{j_1+...+j_n}P_{j_1...j_n}y_1...\hat{y}_{j_1}...\hat{y}_{j_n}...y_m,\,\,\,
\end{equation}

is a cycle representing the fundamental class of the
cohomology algebra.


\section{The spectral sequence}

  In what follows, we give the expression for $\delta$ in the case where k=3.

As mentioned in the introdction, our filtration is one of filterd differentail graded algebras, hence
in this case the identification (\ref{ssk}) becomes :
$$ E_0^{p,q} = \big(\Lambda^{2p}V\oplus\Lambda^{2p+1}V\big)^{p+q} $$
with the product given by: $$(u,v)\otimes(u',v')=(uu',uv'+vu'),\;\;\;  \forall (u,v)\in E_0^{p,q}, \forall (u',v')\in E_0^{p',q'}.$$

On the other hand, since $d_1=d_2=0$ the diffential on $ E_0$ is zero , hence $E_1^{p,q} = E_0^{p,q}$ and so
 the identification obove gives the following diagram

$$\begin{tikzpicture}
  \matrix (m) [matrix of math nodes,row sep=3em,column sep=8em,minimum width=2em]
  {
     E_1^{p,q} & \big(\Lambda^{2p}V\oplus\Lambda^{2p+1}V\big)^{p+q} \\
     E_1^{p+1,q} & \big(\Lambda^{2(p+1)}V\oplus\Lambda^{2(p+1)+1}V\big)^{p+q+1} \\};
  \path[-stealth]
    (m-1-1) edge node [left] {$\delta$} (m-2-1)
            edge  node [above] {$\cong$} (m-1-2)
    (m-2-1.east|-m-2-2) edge
            node [above] {$\cong$} (m-2-2)
    (m-1-2) edge node [left] {$\delta$} (m-2-2)
            ;
    \draw[color=red][<-,dashed] (0.5,-0.5) -- (1,0.5) node[pos=0.5,below] {$d_3$};
    \draw[color=red][<-,dashed] (1.7,-0.5) -- (1.2,0.5) node[pos=0.5,below] {$d_4$};
     \draw[color=red][<-,dashed] (1.9,-0.5) -- (2.4,0.5) node[pos=0.5,below] {$d_3$};
\end{tikzpicture}$$
with $\delta$  defined as follows,
$$\delta(u,v)=(d_3u,d_3v+d_4u)$$

 Let $E_1^p=E_1^{p,\ast}=\bigoplus\limits_{\underset{}{q\geq0}}E_1^{p,q}$ and $E_1^{\ast}=\bigoplus\limits_{\underset{}{p\geq0}}E_1^{p,\ast}$.
 This gives a commutative differential graded algebra $(E_1^{\ast}, \delta)$ wich is the first term of our spectral sequence:
 $$E_2^{p,q} = H^{p,q}(\Lambda V,\delta)\Rightarrow H^{p+q}(\Lambda V,d).$$


\section{Proof of the theorem 1.1}

   Recall that we restrict ourself to the case $k=3$. The approch used here is inspered by that used in \cite{murillo02}.
   Note also that  the subsequent notations imposed us to replace certain somes by pairs and vice-versa. \\

\textbf{For the first inequality}\\

We note first that since by hypothesis, $dimH^N(\Lambda V,d)=1$, the class $\alpha \in E_2^{*,*}$ must survive to $E_{\infty}$.

In what follow we put :
$r = sup\{ k\geq 0, \; \alpha = [\omega _0]  \; with \;  \omega _0 \in \Lambda^{\geq k} V\}.$

  Let then $\omega_0\in \Lambda^{\geq r}V$. We may suppose that $r=2p$ is even (inded, if $r = 2p+1$ is odd, it suffice to rewrite $\omega _0$ with the coordinate in $\Lambda ^{2p}V$ being $0$). More explicily $\omega_0\in(\Lambda^{2p} V\oplus\Lambda^{2p+1} V)\oplus(\Lambda^{2p+2} V\oplus\Lambda^{2p+3} V)\oplus...$,

Since $\mid \omega _0 \mid = N$,  there is  an integer $l$ such that:
 $$ \omega_0 = \omega^0_0+\omega^1_0+...+\omega^l_0 \,\,\,\,\,\,
 with \,\,\,\,\,\,\omega^i_0=(\omega^{i,1}_0,\omega^{i,2}_0)\in\Lambda^{2(p+i)} V\oplus\Lambda^{2(p+i)+1} V$$
We hace successivly:
$$\delta(\omega^i_0)=\delta(\omega^{i,1}_0,\omega^{i,2}_0)=(d_3\omega^{i,1}_0,d_3\omega^{i,2}_0+d_4\omega^{i,1}_0)$$
$$\delta(\omega_0)=\sum_{i=0}^l\delta(\omega^{i,1}_0,\omega^{i,2}_0)=\sum_{i=0}^l(d_3\omega^{i,1}_0,d_3\omega^{i,2}_0+d_4\omega^{i,1}_0)$$
Also, we have $ d\omega_0=d\omega_0^0+d\omega_0^1+ ... +d\omega_0^l$, with:

$$d\omega_0^0=d(\omega_0^{0,1},\omega_0^{0,2})=(d_3\omega_0^{0,1},d_3\omega_0^{0,2}+d_4\omega_0^{0,1})+ ...\in(\Lambda^{2p+2} V\oplus\Lambda^{2p+3} V)\oplus... $$
$$d\omega_0^1=d(\omega_0^{1,1},\omega_0^{1,2})=(d_3\omega_0^{1,1},d_3\omega_0^{1,2}+d_4\omega_0^{1,1})+ ...\in(\Lambda^{2p+4} V\oplus\Lambda^{2p+5} V)\oplus... $$
$$........$$
$$d\omega_0^i=d(\omega_0^{i,1},\omega_0^{i,2})=(d_3\omega_0^{i,1},d_3\omega_0^{i,2}+d_4\omega_0^{i,1})+ ...\in(\Lambda^{2p+2i} V\oplus\Lambda^{2p+2i+1} V)\oplus... $$
Therfore
\begin{align*}
d\omega_0&=(d_3(\omega_0^{0,1}+\omega_0^{1,1}+\omega_0^{2,1}+...)+d_4\omega_0^{0,2}+d_5\omega_0^{0,1}+...,d_3(\omega_0^{0,2}+\omega_0^{1,2}+...)\\
&+d_4(\omega_0^{0,1}+\omega_0^{1,1}+\omega_0^{2,1}+...)+d_5\omega_0^{0,2}+d_6\omega_0^{0,1}+...)
\end{align*}
that is:
$d\omega_0=\delta(\omega_0)+(d_4\omega_0^{0,2}+d_5\omega_0^{0,1}+...,d_5\omega_0^{0,2}+d_6\omega_0^{0,1}+...)$. As  $\delta(\omega_0)=0$\\
we can rewrite:
$$d\omega_0=a_2^0 + a_3^0 + ... + a_{t+l}^0 \,\,\,\,\, with \,\,\,\,\,a_i^0=(a_i^{0,1},a_i^{0,2})\in\Lambda^{2(p+i)} V\oplus\Lambda^{2(p+i)+1} V$$

 Note also that t is a fix integer.
Indeed the degree of $a_{t+l}^0$ is greater or equal than $2(2(p+t+l)+1)$
and it coincides with N + 1, being N the formal dimension.\\
Then $$N+1\geq 2(2(p+t+l)+1)$$
Hence $$t\leq \frac{1}{4}(N-4p-4l-1).$$

In what follows, we take $t$ the largest integer satisfying this enequality.

Now, we have:
\begin{align*}
d^2\omega_0&= da_2^0+da_3^0+ ... +da_{t+l}^0\\
&=  d(a_2^{0,1},a_2^{0,2})+d(a_3^{0,1},a_3^{0,2})+...+d(a_{t+l}^{0,1},a_{t+l}^{0,2})
\end{align*}
with;
\begin{align*}
d(a_2^{0,1},a_2^{0,2})&= d_3(a_2^{0,1},a_2^{0,2})+d_4(a_2^{0,1},a_2^{0,2})+d_5(a_2^{0,1},a_2^{0,2})+...\\
&=  (d_3a_2^{0,1},d_3a_2^{0,2}+d_4a_2^{0,1})+(d_5a_2^{0,1}+d_4a_2^{0,2},d_6a_2^{0,1}+d_5a_2^{0,2})+...
\end{align*}
\begin{align*}
d(a_3^{0,1},a_3^{0,2})&= d_3(a_3^{0,1},a_3^{0,2})+d_4(a_3^{0,1},a_3^{0,2})+d_5(a_23^{0,1},a_3^{0,2})+...\\
&=  (d_3a_3^{0,1},d_3a_3^{0,2}+d_4a_3^{0,1})+(d_5a_3^{0,1}+d_4a_3^{0,2},d_6a_3^{0,1}+d_5a_3^{0,2})+...
\end{align*}
$$........$$

It follows that:
$$d^2\omega_0= (d_3a_2^{0,1},d_3a_2^{0,2}+d_4a_2^{0,1})+(d_5a_2^{0,1}+d_4a_2^{0,2}+d_3a_3^{0,1},d_6a_2^{0,1}+d_5a_2^{0,2}+d_4a_3^{0,1}+d_3a_3^{0,2})+...$$

Since $d^2\omega_0=0$, we have  $(d_3a_2^{0,1},d_3a_2^{0,2}+d_4a_2^{0,1})=\delta(a_2^0)=0$  with
$a_2^0\in \Lambda^{2(p+2)} V\oplus\Lambda^{2(p+2)+1} V $.
Hence $a_2^0$ is a $\delta$-boundary, i.e., there is $b_2\in \Lambda^{2(p+2)-2} V\oplus\Lambda^{2(p+2)-1} V $ such that $a_2^0=\delta(b_2)$. Otherwise the cocycle will not survive to $E_3$ and a fortiori to $E_{\infty}$. 

Consider $\omega_1=\omega_0-b_2$ and reconsider the previous calculation:
\begin{align*}
 d\omega_1&=d\omega_0-db_2\\
 &=(a_2^0+a_3^0+ ... +a_{t+l}^0)-(d_3b_2+d_4b_2+...+d_{t+3}b_2)
\end{align*}
With\\
$d_3b_2=d_3(b_2^1,b_2^2)=(d_3b_2^1,d_3b_2^2)\in \Lambda^{2p+4} V\oplus\Lambda^{2p+5} V$ \\
$d_4b_2=d_4(b_2^1,b_2^2)=(d_4b_2^1,d_4b_2^2)\in \Lambda^{2p+5} V\oplus\Lambda^{2p+6} V $\\
$ .............$

This imply that:
\begin{align*}
 d\omega_1&=a_2^0+a_3^0+ ... +a_{r+l}^0-(d_3b_2^1,d_3b_2^2+d_4b_2^1)+...\\
          &=a_2^0-\delta b_2+a_3^0+ ... +a_{r+l}^0-(d_5b_2^1+d_4b_2^2,d_5b_2^2+...)-...\\
          &=a_3^0-(d_5b_2^1+d_4b_2^2,d_5b_2^2+...)+...
\end{align*}
and then:
$$d\omega_1=a_3^1+a_4^1+...+a_{t+l}^1,\,\,\,\,\, with \,\,\,\,\, a_i^1\in \Lambda^{2(p+i)} V\oplus\Lambda^{2(p+i)+1} V $$
So,
\begin{align*}
d^2\omega_1&=da_3^1+da_4^1+...+da_{t+l}^1\\
           &=d(a_3^{1,1},a_3^{1,2})+d(a_4^{1,1},a_4^{1,2})+...+d(a_{t+l}^{1,1},a_{t+l}^{1,2})\\
           &=(d_3a_3^{1,1},d_3a_3^{1,2}+d_4a_3^{1,1})+(d_5a_3^{1,1}+d_4a_3^{1,2}+d_3a_4^{1,1},d_5a_3^{1,2}+...)+...
\end{align*}

Since $d^2\omega_1=0$, by wordlength reasons, $(d_3a_3^{1,1},d_3a_3^{1,2}+d_4a_3^{1,1})=\delta(a_3^1)=0$. Hence ( for the same reason as before ) $a_3^1$ is a $\delta$-boundary, i.e., there is $b_3\in \Lambda^{2(p+3)-2} V\oplus\Lambda^{2(p+3)-1} V $  such that $\delta(b_3)=a_3^1.$\\
Consider $\omega_2=\omega_1-b_3.$\\
By the same way we show that
$$d\omega_2=a_4^2+a_5^2+...+a_{t+l}^2,\,\,\,\,\, with \,\,\,\,\,a_i^2\in \Lambda^{2(p+i)} V\oplus\Lambda^{2(p+i)+1} V. $$

We continue this process defining inductively $\omega_j=\omega_{j-1}-b_{j+1}$, $j< r+l$  such that:
$$d\omega_j=a_{j+2}^j+a_{j+3}^j+...+a_{t+l}^j,\,\,\,\,\, with \,\,\,\,\,a_i^j\in \Lambda^{2(p+i)} V\oplus\Lambda^{2(p+i)+1} V $$
Also, we have:
$$\omega_{t+l-2}=\omega_{t+l-3}-b_{t+l-1},\,\,\,\,\, with\,\,\,\,\,b_{t+l-1}\in \Lambda^{2(p+t+l-1)-2} V\oplus\Lambda^{2(p+t+l-1)-1} V $$
$$d\omega_{t+l-2}=a_{t+l}^{t+l-2}=\delta(b_{t+l-1})\in \Lambda^{2(p+t+l)} V\oplus\Lambda^{2(p+t+l)+1} V$$
$$d^2\omega_{t+l-2}=da_{t+l}^{t+l-2}=(d_3a_{t+l}^{t+l-2,1},d_3a_{t+l}^{t+l-2,2}+d_4a_{t+l}^{t+l-2,1})+...$$
Since $d^2\omega_{t+l-2}=0$, by wordlength reasons,
$$(d_3a_{t+l}^{t+l-2,1},d_3a_{t+l}^{t+l-2,2}+d_4a_{t+l}^{t+l-2,1})=\delta(a_{t+l}^{t+l-2})=0$$
Hence $a_{t+l}^{t+l-2}$ is a $\delta$-boundary, i.e., there is $b_{t+l}\in \Lambda^{2(p+t+l)-2} V\oplus\Lambda^{2(p+t+l)-1} V $ such that $\delta(b_{t+l})=a_{t+l}^{t+l-2}$.\\
Consider $\omega_{t+l-1}=\omega_{t+l-2}-b_{t+l}$.\\

Note that $|d\omega_{t+l-1}|=|d\omega_{t+l-2}|=N+1$,
but by the hypothesis on $t$, we have:
$$|d(\omega_{t+l-2}-b_{t+l})|=|a_{t+l}^{t+l-2}-\delta(b_{t+l})-(d-\delta)b_{t+l}|=|-(d-\delta)b_{t+l}|>N+1,$$
then $d\omega_{t+l-1}=0$
and so $\omega_{t+l-1}$ can't be a d-boundary.
Indeed suppose that $\omega_{t+l-1}=(\omega_0^0+\omega_0^1+...+\omega_0^l)-(b_2+b_3+...+b_{t+l})$ were a d-boundary,
By wordlength reasons, $\omega_0^0$ would be a $\delta$-boundary, i.e., there is $x\in \Lambda^{2p-2} V\oplus\Lambda^{2p-1} V$ such that $\delta(x)=\omega_0^0$.
Then
    $$\omega_0=\delta(x)+\omega_0^1+...+\omega_0^l$$
Since $\delta(\omega_0)=0$ , we would have $\delta(\omega_0^1+...+\omega_0^l)=0$, but $\omega_0^1+...+\omega_0^l$ is not a $\delta$-boundary.\\
Thus $\omega_{t+l-1}$ is a non trivial cocycle of degree $N$, the formal dimension,  and therefore it represents the fundamental class.\\

Finaly, since $\omega_{t+l-1}\in \Lambda^{\geq r} V$ we have;

$$e_0(\Lambda V,d)\geq r$$

\textbf{For the second inequality}\\

Denote $s=e_0(\Lambda V,d)$ and let $\omega\in\Lambda^{\geq s} V$ be a cocycle representing the generating class $\alpha$  of $H^*(\Lambda V,d)$ . Write
$ \omega=\omega_0+\omega_1+...+\omega_t, \,\, \,\,\,\,\omega_i\in\Lambda^{s+i} V$. We deduce that:
\begin{align*}
d\omega&=(d_3\omega_0+d_3\omega_1+...+d_3\omega_r)+(d_4\omega_0+d_4\omega_1+...+d_4\omega_t)+...\\
       &=\delta(\omega_0 , \omega_1)+...
\end{align*}


Since $d\omega=0$, by wordlength reasons, it follows that $\delta(\omega_0 , \omega_1)=0$.\\
 If $(\omega_0 , \omega_1)$ were a $\delta$-boundary, i.e., $(\omega_0 , \omega_1)$=$\delta(x)$, then
 \begin{align*}
 \omega-dx&=(\omega_0 , \omega_1) +...+ \omega_r-(d_3x+d_4x+...)\\
          &=(\omega_0 , \omega_1)-\delta(x)+(\omega_2+\omega_3+...+\omega_t)-...
 \end{align*}

so $\omega-dx\in \Lambda^{\geq s+2} V$ which contradicts the fact $s=e_0(\Lambda V,d)$.\\
Hence $(\omega_0 , \omega_1)$ represents the generating class of $H^N(\Lambda V,\delta)$.\\

Since $(\omega_0 ,\omega_1)\in \Lambda^{\geq s} V$ we will have $s\leq r$

Hence $$e_0(\Lambda V,d)\leq r $$
We conclude that  $$e_0(\Lambda V,d)= r$$


\section{Some examples and remarks}

\begin{ex}

Let  $(\Lambda V,d)$ be the pure model defined by
$V^{even}=<x_2,x_6>$, \\ $V^{odd}=<y_5, y_{15}, y_{23}>$ , $dx_2=dx_6=0$,
$dy_5=x_2^3$, $dy_{15}= x_2^2x_6^2$  and $dy_{23}=x_6^4 $.

Clearly we have $dimH(\Lambda V,d_3) = \infty $ and $dimH(\Lambda V,d) < \infty $. 

We note also that, since  $N = 37$ is odd, then any representative of the fundamental  class of $(\Lambda V, d)$ will be of the form: $n_1x_2^kx_6^ly_5 + n_2x_2^{k'}x_6^{l'}y_{15} +  n_3x_2^{k''}x_6^{l''}y_{23}$, with $n_1, \; n_2 \; \hbox{and}\;  n_3 \in \mathbb{N}$.

 Using A. Murillo's algorithm (cf. \S 2) the matrix determining the fundamental class  is:
$$A= \begin{pmatrix}
\begin{tikzpicture}
\node (a) at (0,0) {$x_2^2$}; \node (b) at (1,0) {0}; \node
(c)at(0,-0.5){$x_2x_6^2$};

 \node (d) at (1,-0.5) {0};
\node (e) at (0,-1) {0}; \node (f) at (1,-1) {$x_6^3$};

\end{tikzpicture}
\end{pmatrix}$$

So  $\omega _0 =
-x_2^2x^{3}_6y_{15} + x_2x_6^5y_5 \in \Lambda^{\geq6} V$ is an 
generator of this fundamental cohomology class. As in the first example, it is straightforward to verify that there is only two representatives, with $\omega _1 =  x_2^4x_6y_{23} - x_2^2x_6^3y_{15}$ being  the second one. It follow that $e_0(\Lambda V,d) = 6$. 

Remark also that for this model,  $(\omega _0^0, \omega _0^1)=(-x_2^2x^{3}_6y_{15} , x_2x_6^5y_5)\in \Lambda ^6V \oplus \Lambda ^7V$ is a $\delta$-cocycle and in fact $[(\omega _0^0, \omega _0^1)]\in H^N(\Lambda V,\delta)$ is non zero.  $[(x_2x_6^2y_{23} ,0)]$ is another generating class, hence $dimH^N(\Lambda V,\delta) >1$. The algorithm described in remark 1. is applied to $(\omega _0^0, \omega _0^1)$.

On the other hand, $\omega _0$ is  not an  $d_3$-cocycles,   but $0\not =[\omega _1]\in H^N(\Lambda V,d_3)$. Also $0 \not = [x_2x_6^2y_{23}]$ is another generating class of $H^N(\Lambda V,d_3)$, hence $dimH^N(\Lambda V,d_3)>1$.  Application of  the algorithm in the proof of Theorem 5.  in \cite{murillo02} to $\omega _1$ (which is a homogenious $d_3$-cocycle) gives immediatly $\omega _1$ as a representative of the fundamental class of $(\Lambda V,d)$.

Finaly we note also  that $e_0(\Lambda V,d) = 6 \not = (k-2)dimV^{even} + dimV^{odd} =5$.

\end{ex}

\begin{ex}
Let  $(\Lambda V,d)$ be the pure model defined by
$V^{even}=<x_2,x_6>$, \\ $V^{odd}=<y_5, y_{13}, y_{23}>$ , $dx_2=dx_6=0$,
$dy_5=x_2^3$, $dy_{13}= x_2x_6^2$  and $dy_{23}=x_6^4 $.

Clearly we have $dimH(\Lambda V,d_3) = \infty $ and $dimH(\Lambda V,d) < \infty $. 

We note also that, since  $N = 35$ is odd, then any representative of the fundamental  class of $(\Lambda V, d)$ will be of the form: $n_1x_2^kx_6^ly_5 + n_2x_2^{k'}x_6^{l'}y_{13} +  n_3x_2^{k''}x_6^{l''}y_{23}$, with $n_1, \; n_2 \; \hbox{and}\;  n_3 \in \mathbb{N}$.

 Using A. Murillo's algorithm (cf. \S 2) the matrix determining the fundamental class  is:
$$A= \begin{pmatrix}
\begin{tikzpicture}
\node (a) at (0,0) {$x_2^2$}; \node (b) at (1,0) {0}; \node
(c)at(0,-0.5){$x_6^2$};

 \node (d) at (1,-0.5) {0};
\node (e) at (0,-1) {0}; \node (f) at (1,-1) {$x_6^3$};

\end{tikzpicture}
\end{pmatrix}$$

So  $\omega _0 =
-x_2^2x^{3}_6y_{13} + x_6^5y_5 \in \Lambda^{\geq6} V$ is an 
generator of this fundamental cohomology class. Another representative of this class is
$\omega _1 = - x_2^3x_6y_{23} + x_2^2x_6^3y_{13}$.  It is a straightforward calculation to prouve that they are the uniques representatives.   We conclude that $e_0(\Lambda V,d)= 6$.

On the other hand $H^N(\Lambda V, \delta)$ has at least tow generators: $(\omega _0, 0) \in \Lambda ^6 V\oplus \Lambda ^7V$ and $[(0,x_6^2y_{23})]$, hence $dimH^N(\Lambda V, \delta)>1$.  We have also $dimH^N(\Lambda V,d_3) >1$  with  $[\omega _0]$   and  $[x_6^2y_{23}]$ being two generators of $H^N(\Lambda V,d_3)$.   Here the algorithm is applied to $(\omega _0, 0)$ and the one of \cite{murillo02} is applied to $[\omega _0]$.

Note finaly that  $e_0(\Lambda
V,d)=6\neq (k-2)dimV^{even} + dimV^{odd} = 5.$

\end{ex}

\begin{Rq}
It should be noted that the algorithms  that are described in \cite{murillo02} and in Remark 1. are both valid in the previous examples. The previlege of one or the other  depends on $dimH^N(\Lambda V, \delta)$ and $dimH^N(\Lambda V, d_3)$ and also in the expressions of there basis.

On the other hand all the lower bounds for $e_0(\Lambda V,d)$ known up to  now  can be used to relax the application of the algorithm.
\end{Rq}


\end{document}